\documentclass[11pt]{article}
\usepackage{amssymb,amsmath,amsthm,amsfonts}
\usepackage{mathrsfs}
\makeatletter
\def\@fnsymbol#1{\ensuremath{\ifcase#1\or 1\or 2\fi}}
\makeatother

\def\sha{\mathcal{A}}
\def\shb{\mathcal{B}}

\def\shd{\mathcal{D}}

\def\shm{\mathcal{M}}
\def\sho{\mathcal{O}}

\newcommand{\C}{\mathbb{C}}
\newcommand{\N}{\mathbb{N}}
\newcommand{\R}{\mathbb{R}}

\newcommand{\Z}{\mathbb{Z}}
\newtheorem{theorem}{Theorem}[section]
\newtheorem{proposition}[theorem]{Proposition}
\newtheorem{lemma}[theorem]{Lemma}
\newtheorem{corollary}[theorem]{Corollary}

\theoremstyle{definition}
\newtheorem{definition}[theorem]{Definition}

\newtheorem{example}[theorem]{Example}

\newtheorem{remark}[theorem]{Remark}

\begin{document}
\author{ Ana Rita Martins and Teresa Monteiro Fernandes}
\title{Truncated microsupport and hyperbolic inequalities}
\date{\today}

\maketitle \footnote{The research of both authors was supported by
Funda\c c{\~a}o para a Ci{\^e}ncia e Tecnologia and Programa
Ci{\^e}ncia, Tecnologia e Inova\c c{\~a}o do Quadro
Comunit{\'a}rio de Apoio.} 

\begin{abstract}
We prove that the k-truncated microsupport of the specialization
of a complex of sheaves $F$ along a submanifold is contained in
the normal cone to the conormal bundle along the k-truncated
microsupport of $F$. In the complex case, applying our estimates
to $F=\text{R}\cal{H}\text{om}_{\cal{D}}(\cal{M}, \sho)$, where
$\cal{M}$ is a  coherent $\cal{D}$-module, we obtain new estimates
for the truncated microsupport of real analytic and hyperfunction solutions. When
$\cal{M}$ is regular along $Y$ we also obtain estimates for the
truncated microsupport  of the holomorphic solutions of the
induced system along $Y$ as well as for the nearby-cycle sheaf of
$\shm$ when $Y$ is a hypersurface.
\end{abstract}

\section{Introduction and statement of the main results}\label{section:intro}

\hspace*{\parindent} Let $X$ be a real manifold and let $F$ denote an object of the
derived category of abelian sheaves on X.  The microsupport of
$F$, denoted by $SS(F)$, was introduced by M. Kashiwara and P.
Schapira (\cite{K-S1}; \cite{KS2}),  as a subset of the cotangent
bundle $\pi:T^*X\to X$ describing the directions of non
propagation for $F$. The truncated microsupport of a given degree
$k$ (or $k$-truncated microsupport), $SS_k(F)$, defined by  the
same authors, is only concerned by degrees of cohomology up to the
order $k$ and allows us to consider some phenomenon of propagation
in specific degrees. Such notion is particularly useful in the
framework of the theory of linear partial differential equations.
More precisely,  when $F$ is the complex of holomorphic solutions
of a coherent module $\cal{M}$ over the sheaf $\shd_X$ of
holomorphic differential operators on a complex manifold $X$,
$SS_k(F)$ is completely determined as a subset of the
characteristic variety $\text{Char}(\cal{M})$, which itself
coincides with $SS(F)$. In the characteristic case,   interesting
propagation results (cf. \cite{E-K-S}, \cite{T}, \cite{KMS1}) may
be obtained with the truncated microsupport. The truncated
microsupport and its functorial properties were studied in
\cite{KMS1} and  \cite{KMS2}.

It is now natural to study the behaviour of $SS_k(F)$ under
specialization along a submanifold. That is the main purpose of
this work, having in scope the application to $\cal{D}$-modules,
specially to holomorphic solutions of induced systems and to real
analytic solutions.

Let $\mathbf{k}$ be a field. Let $D^b(\mathbf{k}_X)$ denote the
bounded derived category of complexes of sheaves of
$\mathbf{k}$-vector spaces.

Let $M$ be a submanifold of $X$. We shall identify $T_{T^*_M
X}(T^*X)$,  $T^*(T_M X)$ and  $T^*(T^*_M X)$ thanks to the
Hamiltonian isomorphism. Unless otherwise specified, we shall
follow the notations in \cite{K-S1}. In particular, for $F\in
D^b(\mathbf{k}_X)$, $\nu_M(F)$ denotes the specialization of $F$
along $M$, an object of $ D^b(\mathbf{k}_{T_M X}$) and $C_{T^*_M
X}(SS_k(F))$ denotes the normal cone to $T^*_M X$ along $SS_k(F)$.
For a morphism $f:Y\to X$ we shall use $f^{\#}$, a correspondence
which associates conic subsets of $T^*Y$ to conic subsets of
$T^*X$ as well as the operation $\widehat{+}$   which associates
to pairs of  conic closed subsets of $T^*X$   conic closed subsets
of $T^*X$.

The main result of this work is the following:

\begin{theorem}\label{T:26}
Let $M$ be a closed submanifold of $X$ and let $F\in
D^b(\mathbf{k}_X)$. Then: $$SS_k(\nu_M(F))\subset C_{T^*_M
X}(SS_k(F)).$$
\end{theorem}

One difficulty in its proof is that the use of distinguished
triangles is not always convenient because of the shift they
introduce. We also needed to deduce a number of further functorial
properties. Namely, as  an essential step of the proof of this
theorem, we obtain the following estimate:

\begin{theorem}\label{P:6}
Let $Y$ and $X$ be real  manifolds, let $f\colon Y\to X$
be a morphism and let $F\in D^{b}(\mathbf{k}_X)$. Then
$$SS_k(f^{-1} F)\subset f^\#(SS_k(F)).$$
\end{theorem}

Let us denote by $f_d$ and $f_\pi$  the canonical morphisms ($f_d$
was noted by ${}^{t} f'$ in \cite{K-S1}):\begin{center}
$f_\pi:X\times_Y T^*Y\rightarrow T^*Y$ and $f_d:X\times_Y
T^*Y\rightarrow T^*X$.
\end{center}

Regarding $f$ as the composition of a smooth map with a closed
embedding, the proof of Theorem \ref{P:6} relies in two steps. The
first is to apply Proposition 4.4 of \cite{KMS1} which proves the
estimate when $f$ is smooth. The second is  Proposition
\ref{L:30}, where we obtain the estimate $$SS_k(F|_M)\subset
j_dj_\pi^{-1}( SS_k(F)\widehat{+}T_M^*X),$$ when $j$ is  a closed
embedding.

Remark that, in that case, $ j^\#(SS_k(F))= j_dj_\pi^{-1}(
SS_k(F)\widehat{+}T_M^*X).$

In particular, when $f$ is non characteristic with respect to $F$,
we get $$SS_k(f^{-1} F)\subset f_df_\pi^{-1}(SS_k(F)).$$

Namely, when $M$ is non characteristic with respect to $F$, in
other words, $$ SS(F)\cap T_M^*X\subset T^*_M M,$$ we have
$SS_k(F)\widehat{+}T^*_MX=SS_k(F)+T_M^*X$ and
$$j_dj_\pi^{-1}(SS_k(F)+T_M^*X)=j_dj_\pi^{-1}(SS_k(F)).$$

Let now $Y$ be a complex closed smooth hypersurface of a
complex analytic manifold $X$ and assume that $Y$ is defined
as the zero locus of a holomorphic function $f$. Let $\psi_Y$
denote the functor of nearby cycles associated to $Y$. Recall that $Y$
may be regarded as a submanifold $Y'$ of $T_Y X$ by a canonical section
$s$ given by $s$ such that $\psi_Y (F)\simeq s^{-1}\nu_Y (F)$.

Then, Theorem \ref{T:26} entails:

\begin{corollary}\label{P:44}
Let $F\in D^{b}(\mathbf{k}_X)$. Then $$SS_k(\psi_Y (F)) \subset
s_ds_\pi^{-1}(C_{T_Y^*X} (SS_k(F))\widehat{+} T_{Y'}^*(T_Y X)).$$
\end{corollary}

Let us point out that one interesting application of Proposition
\ref{L:30} is the new estimate for the $k$-truncated microsupport
of the tensor product (see Proposition \ref{P:28}).

We end this paper with the application of our results to the
complex $F=R\mathcal{H}om_{\shd_X}(\shm, \sho_X)$ of holomorphic
solutions of a coherent $\shd_X$-module $\shm$ on a complex
manifold $X$(see Section 6.2). Let $\shd_X$ be the sheaf of linear
partial differential operators of finite order and $\sho_X$ the
sheaf of holomorphic functions. Let $Y$ be a complex submanifold
of $X$ and  $j$ be the embedding of $Y$ in $X$. We shall denote by
$\shm_Y$ the induced system, an object of the derived category of
left $\shd_Y$-modules. Recall that, when $\shm$ is regular in the
sense of \cite{KO}, $\shm_Y$ has coherent cohomology. Let
$\tau\colon T_Y X\to Y$ be the projection. Still under the
assumption that $\shm$ is regular along $Y$, one defines a
coherent $\shd_{T_Y X}$-module $\nu_Y(\shm)$, the specialisation
of $\shm$ along $Y$, satisfying a natural isomorphism $$ \nu_Y
(R\mathcal{H}om_{\shd_X}(\shm, \sho_{X}))\simeq
R\mathcal{H}om_{\shd_{T_Y X}}(\nu_Y(\shm), \sho_{T_Y X}).$$
Moreover, if $Y$ has codimension $1$, one defines the nearby-cycle
module $\psi_Y(\shm)$, a coherent $\shd_Y$-module, satisfying a
natural isomorphism $$\psi_Y (F)\simeq
R\mathcal{H}om_{\shd_Y}(\psi_Y(\shm), \sho_Y).$$ We refer
\cite{K5} for the details on these isomorphisms.

Set $V=SS(F)=\text{Char}(\shm)$ and denote by
$V=\bigsqcup_{\alpha} V_{\alpha}$ the (local) decomposition of $V$
in its irreducible components. The notion of orthogonality between
a submanifold $Y$ of $X$ and an involutive subvariety $V$ of
$T^*X$ will be recalled at  Section \ref{S:2}. We recall in Lemma
\ref{L:28} that if $Y, V$ are orthogonal and $V$ is irreducible,
then $V'=j_d(j_\pi^{-1}(V))$ is irreducible and $\pi(V)$ has the
same codimension of $\pi'(V')$. Here, $\pi': T^*Y\to Y$ denotes
the projection.

As a consequence of Theorem \ref{T:26} together with the results
of \cite{K5} we obtain:
\begin{theorem}\label{P:146}
Let $\cal{M}$ be a coherent $\shd_X$-module. Then:
$$SS_k(R\mathcal{H}om_{\tau^{-1}\shd_X}(\tau^{-1}\shm,
\nu_Y(\sho_X)))\subset C_{T_Y^*X}( SS_k(F)).$$ If, moreover,
$\shm$ is regular along $Y$ in the sense of \cite{KO} we have:
$$SS_k(R\mathcal{H}om_{\shd_{T_Y X}}(\nu_Y(\shm), \sho_{T_Y
X}))\subset C_{T_Y^*X}( SS_k(F)).$$
\end{theorem}

From the preceding theorem, the results of \cite{K5} and Corollary
\ref{P:44} we obtain:

\begin{corollary}\label{P:45}
Assume that $\cal{M}$ is regular along $Y$ in the sense of
\cite{KO}. Then $$SS_k(R\mathcal{H}om_{\shd_Y}(\psi_Y(\shm),
\sho_Y))\subset s_ds_\pi^{-1}(C_{T_Y^*X} (SS_k(F))\widehat{+}
T_{Y'}^*(T_Y X)).$$
\end{corollary}

Furthermore, Proposition \ref{L:30} together with the results of \cite{K5} and
Theorem 6.7 of \cite{KMS1} entails:

\begin{theorem}\label{P:46}
Assume that $\cal{M}$ is regular along $Y$ in the sense of
\cite{KO}. Then: $$SS_k(R\mathcal{H}om_{\shd_Y}(\shm_Y,
\sho_Y))\subset j_dj_\pi^{-1}( SS_k(F)\widehat{+}T_Y^*X).$$  If,
moreover, $Y$ is non characteristic for $\shm$ and $Y$ is
orthogonal to each $V_{\alpha}$ such that  codim $\pi
(V_{\alpha})\leq k$, the preceding inclusion becomes an equality,
for every $i\leq k$: $$SS_i(R\mathcal{H}om_{\shd_Y}(\shm_Y,
\sho_Y))= j_dj_\pi^{-1}(SS_i(F)).$$
\end{theorem}

Recall that M. Kashiwara has proven in \cite{K6} that, when $Y$ is
non characteristic for $\shm$,
$$SS(R\mathcal{H}om_{\shd_Y}(\shm_Y, \sho_Y))=j_d
j_\pi^{-1}(SS(F)).$$  The  condition of orthogonality  is required
in  Theorem \ref{P:46} in order to have the analogous equality up
to a given degree $k$.

Let us now assume that the complex manifold $X$ is the
complexified of a real analytic manifold $M$. Denote by  $\sha_ M$
the sheaf of real analytic functions on $M$ and by $j$ the
embedding of $M$ in $X$.

Another important application of Theorem \ref{P:6} is:

\begin{proposition}\label{P:7}
Let $\shm$ be a coherent $\shd_X$-module. Then we have the
estimate: $$SS_k (R\mathcal{H} om_{\shd_X}(\shm, \sha_ M))\subset
j_dj_{\pi}^{-1}(SS_k(F)\widehat{+}T_M^*X).$$
\end{proposition}
Let $\shb_M$ denote the sheaf of Sato's hyperfuntions on $M$.
As an immediate consequence of Proposition \ref{P:7} together with
Theorem 6.7 of \cite{KMS1} we get:

\begin{corollary}\label{P:118}
Let $\shm$ be an coherent $\shd_X$-module. Assume that
$$SS_k(F)\cap T^*_M X\subset T^*_X X.$$ Then, $$\tau^{\leq k}(R\mathcal{H}
om_{\shd_X}(\shm, \sha_ M))\simeq\tau^{\leq k}(R\mathcal{H}
om_{\shd_X}(\shm, \shb_ M)).$$  In particular $$SS_k(R\mathcal{H}
om_{\shd_X}(\shm, \sha_ M))=SS_k(R\mathcal{H}
om_{\shd_X}(\shm, \shb_ M))\subset j_d j_{\pi}^{-1} (SS_k(F))$$
$$\subset j_d j_{\pi}^{-1}((\underset{codim  Y_{\alpha}<    
k}{\bigcup} V_{\alpha}) \cup (\underset{codim Y_{\alpha}=
k}{\bigcup} T_{Y_{\alpha}}^*X)) .$$
\end{corollary}

We shall illustrate this corollary with an example (see Example
\ref{eg}) of a  propagation phenomenon for real analytic solutions
of a class of non elliptic differential operators, which, as far
as we know, is new.

When
$\shm$ is elliptic, in other words, $$ SS(\shm)\cap T^*_M X\subset
T_X^*X$$, we get the
estimate: $$\text{For any k}, \ SS_k (R\mathcal{H}
om_{\shd_X}(\shm, \sha_ M))= SS_k (R\mathcal{H} om_{\shd_X}(\shm,
\shb_ M))$$ $$\subset j_d j_{\pi}^{-1}((\underset{codim
Y_{\alpha}< k}{\bigcup} V_{\alpha}) \cup (\underset{codim
Y_{\alpha}= k}{\bigcup} T_{Y_{\alpha}}^*X)) .$$

We thank M. Kashiwara and P. Schapira  for the useful discussions
through the preparation  of this work.

\section{Notations}

\hspace*{\parindent} We will mainly follow the notations in \cite{K-S1}.

Let $X$ be a real manifold. We denote by
$\tau:TX\rightarrow X$ the tangent bundle to $X$ and by
$\pi:T^*X\rightarrow X$ the cotangent bundle. We identify $X$ with
the zero section of $T^*X$. Given a smooth submanifold $Y$ of $X$,
$T_YX$ denotes the normal bundle to $Y$ and $T_Y^*X$ the conormal
bundle. Given a submanifold $Y$ of $X$ and a subset $S$ of $X$ we
denote by $C_Y(S)$ the normal cone to $S$ along $Y$, a closed
conic subset of $T_YX$.

Let $f:X\rightarrow Y$ be a morphism of manifolds. We denote by
\begin{center}
$f_\pi:X\times_Y T^*Y\rightarrow T^*Y$ and $f_d:X\times_Y
T^*Y\rightarrow T^*X$
\end{center}
the associated morphisms.

Given a subset $A$ of $T^*X$, we denote by $A^a$ the image of $A$
by the antipodal map $$a:(x,\xi)\mapsto (x;-\xi).$$ The closure of
$A$ is denoted by $\overline{A}$. For a cone $\gamma\subset TX$,
the polar cone $\gamma^\circ$ to $\gamma$ is the convex cone in
$T^*X$ defined by
\begin{center}
$\gamma^\circ=\{(x;\xi)\in TX; x\in\pi(\gamma)$ and $\langle
v,\xi\rangle\geq 0$ for any $(x;v)\in\gamma \}.$
\end{center}

Given conic subsets $A$ and $B$ of $T^*X$, the operations $A+B$
and $A\widehat{+}B$ are defined in \cite{K-S1} and will be
recalled in section \ref{rnc}.

Given an open subset $\Omega$ of $X$, as in \cite{K-S1}, we
denote by $N^*(\Omega)$ the conormal cone to $\Omega$.

When $X$ is an open subset of a real finite-dimensional vector
space $E$ and $\gamma$ is a closed convex cone (with vertex at 0)
in $E$, we denote by $X_\gamma$ the open set $X$ endowed with
the induced $\gamma$-topology of $E$.

Let $\mathbf{k}$ be a field. We denote by $D(\mathbf{k}_X)$ the
derived category of complexes of sheaves of $\mathbf{k}$-vector
spaces on $X$ and by $D^b(\mathbf{k}_X)$ the full subcategory of
$D(\mathbf{k}_X)$ consisting of complexes with bounded
cohomologies.

For $k\in\Z$, we denote by $D^{\geq k}(\mathbf{k}_X)$ (resp.
$D^{\leq k}(\mathbf{k}_X)$) the full additive subcategory of
$D^b(\mathbf{k}_X)$ consisting of objects $F$ satisfying
$H^j(F)=0$ for any $j< k$ (resp. $H^j(F)=0$ for any $j> k$). The
category $D^{\geq k+1}(\mathbf{k}_X)$ is  denoted by
$D^{> k}(\mathbf{k}_X)$.

Given an object $F$ of $D^b(\mathbf{k}_X)$ and  a submanifold $M$
of $X$, $\nu_M(F)$ denotes the specialization of $F$ along $M$, an
object of $D^b(\mathbf{k}_{T_M X})$.

Let $F$ be an object of $D^b(\mathbf{k}_X)$; we denote by  $SS(F)$
its microsupport, a closed $\R^+$-conic involutive subset of
$T^*X$. For $p\in T^*X$, $D^b(\mathbf{k}_X;p)$ denotes the
localization of $D^b(\mathbf{k}_X)$ by the full triangulated
subcategory consisting of objects $F$ such that $p\notin SS(F)$.

Let $X$ be a finite-dimensional complex manifold. We denote by
$\mathcal{O}_X$ the sheaf of holomorphic functions, by
$\mathcal{D}_X$ the sheaf  of linear holomorphic differential
operators of finite order and by  $\mathcal{D}_X(\cdot)$ the
filtration by the order. Given  a coherent $\shd_X$-module $\shm$,
we denote by $\text{Char}(\shm)$ its characteristic variety.

Let $Y$ be a closed submanifold, let $\tau$ be the projection of
$T_Y X$ on $Y$ and let $V^{\cdot}_Y$ denote the V-filtration on
$\mathcal{D}_X$ with respect to $Y$.  Let $\shd_{[T_YX]}$ denote
the sheaf of differential operators on $T_Y X$ with polynomial
coefficients with respect to the fibers of $\tau$. Let $\theta$
denote the Euler operator on $T_YX$. Recall that $\shm$ is regular
along $Y$ if for any local section $u$ of $\shm$ there exists a
non trivial polynomial $b$ of degree $m$ such that $$b(\theta)u\in
(V^1_Y(\shd_X)\cap \shd_X (m))u.$$ Following Kashiwara in \cite{K5}, given an appropriate good
$V^{\cdot}_Y$-filtration on $\shm$, the specialized of $\shm$
along $Y$, $\nu_Y(\shm)$,  is the coherent $\shd_{T_Y X}$-module generated by the associated graded module.
When $Y$ is a hypersurface, one defines a coherent $\shd_Y$-module, the nearby-cycles module $\psi_Y(\shm)$, as the degree zero homogeneous term of that graded module.
\section{Review on normal cones in cotangent bundles}\label{rnc}

\hspace*{\parindent} For the reader's convenience  we shall recall here some operations
on conic subsets in cotagent bundles defined on \cite{K-S1}.

Let $X$ be a real manifold, $(x)$ a system of local coordinates on
$X$ and denote by $(x;\xi)$ the associated coordinates on $T^*X$.
Given two conic subsets $A$ and $B$ of $T^*X$, one defines the sum
\begin{center}
$A+B=\{(x;\xi)\in T^*X; \xi=\xi_1+\xi_2$ for some $(x;\xi_1)\in A$
and $(x;\xi_2)\in B\}$.
\end{center}

When $A$ and $B$ are closed,  $A\widehat{+}B$ is the
closed conic set containing $A+B$,  described as follows:
$(x_0;\xi_0)$ belongs to $A\widehat{+}B$ if and only if there
exists sequences $\{(x_n;\xi_n)\}_n$ in $A$ and
$\{(y_n;\eta_n)\}_n$ in $B$ such that: $$\begin{cases} x_n,y_n
\xrightarrow[n]{} x_0,\\ \xi_n+\eta_n \xrightarrow[n]{} \xi_0,\\
|x_n-y_n||\xi_n|\xrightarrow[n]{} 0.
\end{cases}$$

Let  $M$  be a submanifold of $X$. Let $(x',x'')$ be a system of
local coordinates on $X$ such that $M=\{(x',x''); x'=0\}$ and let
$(x',x'';\xi',\xi'')$ denote the associated coordinates on $T^*X$.
Given a subset $\Lambda$ of $T^*X$ we describe the normal cone to
$\Lambda$ along $T_M^*X$, $C_{T_M^*X}(\Lambda)$, as follows:
$(x'_0,x''_0;\xi'_0,\xi''_0)\in C_{T_M^*X}(\Lambda)$ if and only
if there exist sequences of real positive numbers $\{c_n\}_n$ and
$\{(x'_n,x''_n;\xi'_n,\xi_n'')\}_n$ in $\Lambda$ such that:
$$\begin{cases} (x'_n,x''_n;\xi'_n,\xi_n'')\xrightarrow[n]{}
(0,x''_0;\xi'_0,0),\\ c_n(x'_n,\xi''_n)\xrightarrow[n]{}
(x'_0,\xi''_0).
\end{cases}$$

Thanks to the Hamiltonian isomorphism, one gets an embedding of
$T^*M$ into $T_{T_M^*X}T^*X$, and, for a conic subset $\Lambda$ of
$T^*X$, the set $T^*M\cap C_{T_M^*X}(\Lambda)$ is described as
follows: $(x'_0,x''_0;\xi'_0,\xi''_0)\in T^*M\cap
C_{T_M^*X}(\Lambda)$ if and only if there exists a sequence
$\{(x'_n,x''_n;\xi'_n,\xi_n'')\}_n$ in $\Lambda$ such that:
$$\begin{cases} (x''_n;\xi_n'')\xrightarrow[n]{}
(x''_0;\xi''_0),\\ |x'_n|\xrightarrow[n]{} 0\\
|x'_n||\xi'_n|\xrightarrow[n]{} 0.
\end{cases}$$

Let $f :Y\to X$ be a morphism of manifolds. The notion of
$f^{\#}$, a  correspondence introduced in \cite{K-S1} associating
conic subsets of $T^*Y$ to conic subsets of $T^*X$, is rather
complicated and we refer the reader to \cite{K-S1} for the
details. We just recall the following results:

\begin{proposition}\label{P:132}(cf. Proposition 6.2.4 of \cite{K-S1})
Let $\Lambda$ be a conic subset of $T^*X$.

(i) Assume that $f: M\to X$ is a closed embedding. Then, $$f^{\#}
(\Lambda)=T^*M\cap C_{T^*_M X}(\Lambda).$$

(ii) Let $(x)$ (resp. $(y)$) be a system of local coordinates on
$X$ (resp. $Y$) and let $(x;\xi)$ (resp.$(y;\eta)$) be the
associated coordinates on $T^*X$ (resp. $T^*Y$). Then

$(y_0; \eta_0)\in f^{\#}(\Lambda)$ if and only if there exist a
sequence $\{(x_n;\xi_n)\}_n$ in $\Lambda$, a sequence $\{y_n\}_n$
in $Y$ such that $$ y_n\xrightarrow[n]{} y_0, x_n\xrightarrow[n]{}
f(y_0), ({}^t f'(y_n)\cdot\xi_n)\xrightarrow [n]{}\eta_0, \vert
x_n-f(y_n)\vert \vert \xi_n\vert\xrightarrow[n]{} 0. $$
\end{proposition}

We shall also need the following description of $j^{\#}$ when $j$
is an embedding:

\begin{lemma}\label{L:32}
Let  $M$ be a closed submanifold of $X$ and let $j$ denote the
embedding of $M$ in $X$. Let $\Lambda$ be a closed conic subset of
$T^*X$. Then: $$j_dj_\pi^{-1}(\Lambda\widehat{+}T_M^*X)= T^*M\cap
C_{T^*_M X}(\Lambda),$$ where we identify $T_{T^*_M X}T^*X$ and
$T^*T_M X$ by the Hamiltonian isomorphism.
\end{lemma}

\begin{proof}[\textbf{Proof}]
It is enough to prove that
$$j_dj_\pi^{-1}(\Lambda\widehat{+}T_M^*X)=j^{\#}(\Lambda).$$ Let
$p\in j_dj_\pi^{-1}(\Lambda\widehat{+}T_M^*X)$ and let $(x',x'')$
be a system of local coordinates on $X$ in a neighborhood of $p$
such that $M=\{(x); x'=0\}$. Let $(x;\xi)$ denote the associated
coordinates on $T^*X$. Suppose $p=(x_0'';\xi''_0)$.

Then there exists $\xi_0'$ such that $(0,x''_0;\xi_0',\xi''_0)\in
\Lambda\widehat{+}T_M^*X$. By definition of $\widehat{+},$ there
exist sequences $\{(x_n',x''_n;\xi'_n,\xi''_n)\}_n$ in $\Lambda$
and $\{(0,y_n'';\eta_n',0)\}_n$ in $T_M ^*X$ such that
$$\begin{cases}(x_n',x''_n), (0,y_n'')\xrightarrow[n]{}
(0,x''_0),\\ \xi_n''\xrightarrow[n]{} \xi''_0,\\
\xi'_n+\eta'_n\xrightarrow[n]{} \xi'_0,\\
|(x_n',x''_n)-(0,y_n'')||(\xi'_n,\xi''_n)|\xrightarrow[n]{}0.\end{cases}$$
Hence $$\begin{cases}x''_n\xrightarrow[n]{} x''_0,
\\x_n'\xrightarrow[n]{}0,\\ \xi_n''\xrightarrow[n]{} \xi''_0,\\
|x_n'||\xi_n|\xrightarrow[n]{}0
\end{cases}$$ and $(x''_0;\xi''_0)\in j^{\#}(\Lambda) $.

Conversely, let  $p\in j^{\#}(\Lambda)$, $p=(x_0'';\xi''_0)$. Then
there exists a sequence $\{(x'_n,x''_n;\xi'_n,\xi''_n)\}_n$ in
$\Lambda$ such that
$$\begin{cases}(x''_n;\xi''_n)\xrightarrow[n]{} (x''_0;\xi''_0),\\
x_n'\xrightarrow[n]{}0,\\ |x_n'||\xi_n|\xrightarrow[n]{}0
.\end{cases}$$ The sequences $\{(x'_n,x''_n;\xi'_n,\xi''_n)\}_n$
in $\Lambda$ and $\{(0,x''_n;-\xi'_n,0)\}$ in $T_M^*X$ satisfy the
necessary conditions  so that  $(0,x_0'';0,\xi''_0)\in
\Lambda\widehat{+}T_M^*X$, hence $(x''_0;\xi''_0)\in
j_dj_\pi^{-1}(\Lambda\widehat{+}T_M^*X)$.
\end{proof}

\begin{lemma}\label{L:133}
Let $\Lambda$ be a closed conic subset of $T^*X$ and $M$ a closed
submanifold of $X$. One has:
$$(\Lambda\widehat{+}T_M^*X)\widehat{+}T_M^*X=\Lambda\widehat{+}T_M^*X.$$
\end{lemma}

\begin{proof}[\textbf{Proof}]
Let $(x',x'')$ be a system of local coordinates on $X$ such that
$M=\{(x',x'');x'=0\}$ and let $(x',x'';\xi',\xi'')$ be the
associated coordinates on $T^*X$.

Let $(x_0,;\xi_0)\in(\Lambda\widehat{+}T_M^*X)\widehat{+}T_M^*X$,
then there exists sequences $\{(x_n;\xi_n)\}_n$ and
$\{(y_n;\eta_n)\}_n$ in $\Lambda\widehat{+}T_M^*X$ and $T_M^*X$,
respectively, such that
$$\begin{cases}x_n,y_n\xrightarrow[n]{}x_0,
\\ \xi_n+\eta_n\xrightarrow[n]{}\xi_0,\\ |x_n-y_n||\xi_n|\xrightarrow[n]{} 0.\end{cases}$$

For each $n\in \mathbb{N}$, since $(x_n;\xi_n)\in
\Lambda\widehat{+}T_M^*X$ , there exist sequences
$\{(x_m^n;\xi_m^n)\}_m$ in $\Lambda$ and $\{(y_m^n;\eta_m^n)\}_m$
in $T_M^*X$ such that $$\begin{cases}x_m^n,
y_m^n\xrightarrow[m]{}x_n,
\\ \xi_m^n+\eta_m^n\xrightarrow[m]{} \xi_n,\\
|x_m^n-y_m^n||\xi_m^n|\xrightarrow[m]{}0.\end{cases}$$

Hence we  can find subsequences $\{(x_k;\xi_k)\}_k$ and
$\{(y_k;\eta_k))  \}_k$  of $\{(x_m^n;\xi_m^n)\}_{n,m}$ and
$\{(y_m^n;\eta_m^n+\eta_n)\}_{m,n}$, respectively, such that
$$\begin{cases}x_k, y_k\xrightarrow[k]{} x_0, \\
\xi_k+\eta_k\xrightarrow[k]{}\xi_0,\\
|x_k-y_k||\xi_k|\xrightarrow[k]{}0,\end{cases}$$ which gives
$(x_0;\xi_0)\in \Lambda\widehat{+} T_M^*X$.

Conversely, since $\pi(\Lambda\widehat{+}T_M^*X)\subset M$,
$(\Lambda\widehat{+}T_M^*X)\subset(\Lambda\widehat{+}T_M^*X)+T_M^*X\subset
(\Lambda\widehat{+}T_M^*X)\widehat{+}T_M^*X.$
\end{proof}

Let us now  assume that $X$ is an open subset of $\mathbb{R}^{n}$
with the coordinates $(x)=(x_1,...,x_n)$  and that $M$ is the
submanifold  $\{(x',x'')\in X; (x')=(x_1,...,x_d)=0\}$. Let
$\delta>0$ and let $\gamma$ be the closed convex proper cone given
by: $$\gamma=\{(x', x''); x_n\leq -\frac{1}{\delta} \vert
(x',x_{d+1},...x_{n-1}) \vert\}.$$

Hence $$\gamma^{\circ a}=\{(\xi',\xi''); \xi_n\geq \delta
\vert (\xi',\xi_{d+1},...\xi_{n-1}) \vert \}.$$

Therefore  $(x+\gamma)\cap M=x+(\gamma\cap M)$, for each $x\in M$.
Let $\mathbb{R^+}$ denote the set of real positive numbers and let
us introduce the following notation: for any $\lambda\in
\mathbb{R^+}$ $$ \gamma_{\lambda}= \{(x',x')\in X; (\lambda^{-1}
x',x'')\in \gamma\} $$ $$V_{\lambda}=\{(x',x''); (\lambda^{-1}
x',x'')\in V\}.$$

Remark that if  $\lambda<1$, $\text{Int}(\gamma_{\lambda} ^{\circ
a})\supset \gamma ^{\circ a}$.

\begin{lemma}\label{L:131}
Let $\Lambda$ be a conic closed subset of $T^*X$.

Let $x''\in M\cap\pi(\Lambda)$ and assume that there is a compact
neighborhood $V$ of $x''$ such that $$(V\times \gamma^{\circ
a})\cap (\Lambda\widehat{+}T_M^*X) \subset T^*_{X} X.$$ Then,
there exists a real positive number $C$ such that for any
$\lambda$ and $\epsilon$ satisfying $0<\lambda, \epsilon<C,$
$$(V_{\lambda{} \epsilon} \times\gamma^{\circ a} _{\lambda})\cap
\Lambda\subset T^*_X X.$$
\end{lemma}

\begin{proof}[\textbf{Proof}]
We shall argue by contradiction. Therefore, we can find  sequences
$(\lambda_l)_{l\in \mathbb{N}}, (\epsilon_l)_{l\in \mathbb{N}}$ of
positive numbers converging to $0$, $(x'_l ,x''_l;
\xi'_l,\xi''_{l})_{ t\in  \mathbb{N}}$ in $\Lambda$,
$(\xi'_l,\xi''_{l})\neq (0,0)$,  such that $\vert x'_{l}\vert\leq
{C' \epsilon_{l}\lambda_{l}}$, for some positive constant $C'$
only depending on $V$, and $(0, x''_{l};
\lambda_{l}\xi'_{l},\xi''_{l})\in V\times \gamma^{\circ a}$.

Since the $n$-component $(\xi_ l)_{l\in \mathbb{N}}$ is positive,
after dividing  $(\xi'_l,\xi''_{l})$ by $\xi_{l,n}$, we may assume
that $\xi_{l,n}=1$ and that $(\lambda_{l}\xi'_{l}, \xi''_{l})$ is
a bounded sequence. Since $\vert\xi'_{l}\vert \vert x'_l \vert\leq
C' \epsilon_{l}\lambda_{l}\vert\xi'_{l}\vert $ we get that
$(\vert\xi'_{l}\vert \vert x'_l \vert)_{l} $ converges to $0$.
Moreover, since $x''_{l}$ is bounded, we may assume that $x''_{l}$
converges to some $\check{x''}\in V\cap M$ and that
$(\lambda_{l}\xi'_{l}, \xi''_{l})$ converges to some $({\xi'_{0}},
{\xi''_{0}})\in \gamma^{\circ a}$, with $\xi_{0n}=1$. Considering
the sequences $(x'_l ,x''_l ;\xi'_l , \xi''_{l})_{l\in
\mathbb{N}}\in\Lambda$ and  $(0, x''_{l};
-\xi'_{l}+\lambda_{l}\xi'_{l},0)\in T^*_M X$ we get that
$(0,\check{x''};{\xi'_{0}}, {\xi''_{0}})\in (V\times \gamma^{\circ
a})\cap(\Lambda\widehat{+}T_M^*X)$, which entails $\xi_{0n}=0$, a
contradiction.
\end{proof}

Let $\Omega$ be an open subset of $X$. We shall now recall the
notion of conormal cone to $\Omega$, $N^*(\Omega)$. It is the
subset of $T^*X$ defined as follows:

Given $x\in X$, we denote by $N_x(\Omega)$ the subset of $T_x X$
consisting of vectors $v\neq 0$ such that, in a local chart in a
neighborhood of $x$, there exist an open cone $\gamma$  containing
$v$ and a neighborhood $U$ of $x$ such that $$U\cap((\Omega\cap
U)+\gamma)\subset \Omega.$$ Note that, in particular,
$N_x(\Omega)=T_x X$ if and only if $x\notin \overline{\Omega}$ or
$x\in \Omega$. We denote by $N(\Omega)$ the open convex cone of
$TX$:  $$N(\Omega)=\bigcup_{x\in X}N_x(\Omega),$$ and call it the
strict normal cone to $\Omega$.

Finally  $N^*(\Omega)$, the conormal cone to $\Omega$, is given by
$$N^*(\Omega)=\bigcup_{x\in X}(N^*_x(\Omega)),$$ where, for
each $x\in \Omega$, $N^*_x(\Omega)=(N_x(\Omega))^\circ$.

\section{Review on the truncated microsupport}

\hspace*{\parindent} We shall now recall equivalent definitions of the truncated
microsupport, following ~\cite{KMS1}.

Given $(x_{0},\xi_{0})\in \mathbb{R}^n\times(\mathbb{R}^n)^*$ and
$\varepsilon\in\mathbb{R}$ we set:
$$H_{\varepsilon}(x_{0},\xi_{0})=\{x\in\mathbb{R}^n;\langle
x-x_{0},\xi_{0} \rangle> -\varepsilon\},$$ and if there is no risk
of confusion we will write $H_{\varepsilon}$ instead of
$H_{\varepsilon}(x_{0},\xi_{0})$.

\begin{proposition}\label{PP:1}
Let $X$ be a real analytic manifold and let $p\in T^*X$. Let $F\in
D^b(\mathbf{k}_X)$, $k\in\mathbb{Z}$ and $\alpha\in\mathbb{N}\cup
\{ \infty, \omega\}$. Then the following conditions are
equivalent:

$(i)_{k}$ There exist $F'\in D^{>k}(\mathbf{k}_X)$ and an
isomorphism $F\simeq F'$ in $D^b(\mathbf{k}_X;p)$;

$(ii)_{k}$ There exist $F'\in D^{>k}(\mathbf{k}_X)$ and a morphism
$F'\rightarrow F$ in $D^b(\mathbf{k}_X)$ which is an isomorphism
in $D^b(\mathbf{k}_X;p)$;

$(iii)_{k,\alpha}$ There exists an open conic neighborhood $U$ of
$p$ such that for any $x\in \pi(U)$ and any $\mathbb{R}$-valued
$C^\alpha$-function $\varphi$ defined on a neighborhood of $x$
such that $\varphi(x)=0$, $d\varphi(x)\in U$, one has
\begin{center}
$H^j_{\{\varphi \geq 0\}}(F)_x=0,$ for any $j\leq k.$
\end{center}

When $X$ is an open subset of $\mathbb{R}^n$ and
$p=(x_{0},\xi_{0})$, the above conditions are also equivalent to:

$(iv)_{k}$ There exist a proper closed convex cone
$\gamma\subset\mathbb{R}^n$, $\varepsilon>0$ and an open
neighborhood $W$ of $x_{0}$ with $\xi_{0}\in Int(\gamma^\circ)$
such that $(W+\gamma^a)\cap \overline{H_{\varepsilon}}\subset X$
and
\begin{center}
$H^j(X; \mathbf{k}_{(x+\gamma^a)\cap H_{\varepsilon}}\otimes
F)=0,$ for any $j\leq k, x\in W.$
\end{center}
\end{proposition}

\begin{remark}\label{R:3}
Note that when $X$ is an open subset of $\mathbb{R}^n$ and
$p=(x_{0},\xi_{0})$, the equivalent conditions of Proposition
\ref{PP:1} are also equivalent to:

There exists some $F'\in D^b(\mathbf{k}_X)$ isomorphic to $F$ in a
neighborhood of $x_0$ and a closed proper convex cone  $\gamma$ in
$E$, with $0\in \gamma$ and $\xi_0\in Int\gamma^{\circ a}$, such
that $R\phi_{\gamma *}( F')\in D^{> k}(\mathbf{k}_{X_\gamma})$.
\end{remark}

\begin{definition}
Let $F\in D^b(\mathbf{k}_X)$. We define the closed conic subset
$SS_k(F)$ of $T^*X$ by: $p\notin SS_k(F)$ if and only if $F$
satisfies the equivalent conditions in the preceding Proposition.
\end{definition}

We shall need  the following properties of the truncaded
microsupport  also proved in \cite{KMS1}:

(i) Given a distinguished triangle $F'\rightarrow F \rightarrow
F''\xrightarrow[]{+1}$, one has
\begin{equation}
SS_k(F)\subset SS_k(F')\cup SS_k(F''),
\end{equation}

\begin{equation}
(SS_k(F')\backslash SS_{k-1}(F''))\cup (SS_k(F'')\backslash
SS_{k+1}(F'))\subset SS_k(F).
\end{equation}

(ii) For any $F\in D^b(\mathbf{k}_X)$, one has
\begin{equation}\label{tau}
SS_k(F)\cap T_X^*X=\pi(SS_k(F))=supp(\tau^{\leq k}(F)).
\end{equation}

\begin{proposition}\label{P:21}
Let $X$ and $Y$ be two manifolds. Then for $F\in
D^b(\mathbf{k}_X)$, $G\in D^b(\mathbf{k}_Y)$ and $k\in\mathbb{Z}$,
one has: $$SS_k(F\boxtimes G)=\bigcup_{i+j=k}SS_i(F)\times
SS_j(G).$$
\end{proposition}

\begin{proposition}\label{P:2}
Let $Y$ and $X$ be two manifolds, let $f\colon Y\to X$
be a morphism and let $G\in D^{b}(\mathbf{k}_Y)$ such that $f$ is
proper on the support of $G$. Then for any $k\in\N$,

\begin{equation}\label{E:1}
SS_k (Rf_* (G))\subset f_{\pi}f_d^{-1}(SS_k(G)).
\end{equation}
The equality holds in the case f is a closed embedding.
\end{proposition}

\begin{proposition}\label{P:3}
Let $Y$ and $X$ be two manifolds and let $f\colon Y\to
X$ be a smooth morphism. Let $F\in D^b (\mathbf{k}_X)$. Then

\begin{equation}\label{E:8}
SS_k (f^{-1}F)= f_df_{\pi}^{-1}(SS_k (F)).
\end{equation}

\end{proposition}
To end this section, we shall prove the following
characterizations of the truncated microsupport not included in
\cite{KMS1}, which will be useful in the sequel.

\begin{lemma}\label{L:13}
Let $E$ be a real finite-dimensional vector space, $X$ an open
subset of $E$ and let $F\in  D^b (\mathbf{k}_X)$. Let $U$ be an
open subset of $X$ and $\gamma$ be a closed convex proper cone in
$E$ with $0\in \gamma$. Assume that $$SS_k(F)\cap(U\times
\text{Int} (\gamma^{\circ a}))=\emptyset.$$

Then, given $(x_0,\xi_0)\in U\times Int (\gamma^{\circ a})$,
$\varepsilon> 0$ and an open subset $W\subset X$ such that
$(W+\gamma)\cap H_{\varepsilon}\subset\subset U$, one has:
\begin{equation}\label{E:14}
H^j(X; \mathbf{k}_{(x+\gamma)\cap H_{\varepsilon}}\otimes F)=0,\
for\ any \ x\in W+\gamma \ and \ j\leq k.
\end{equation}

\end{lemma}

\begin{proof}[\textbf{Proof}]
We may assume that $X$ is an open subset of $\R^n$.

Let $(x_0;\xi_0)\in U\times Int (\gamma^{\circ a})$, $\varepsilon>
0$ and $W\subset X$  be an open subset such that $(W+\gamma)\cap
H_{\varepsilon}\subset\subset U$. Let us prove (\ref{E:14}).

By the microlocal cut-off lemma (Proposition 5.2.3 of
~\cite{K-S1}), we have a distinguished triangle
$$\phi_\gamma^{-1}R\phi_{\gamma *}F\rightarrow F\rightarrow
G\xrightarrow{+1},$$ with $SS(G)\cap(X\times \text{Int}
(\gamma^{\circ a}))=\emptyset$. Therefore, setting
$F'=\phi_\gamma^{-1}R\phi_{\gamma *}F$, one has $$H^j(X;
\mathbf{k}_{(x+\gamma)\cap H_{\varepsilon}}\otimes F)\simeq H^j(X;
\mathbf{k}_{(x+\gamma)\cap H_{\varepsilon}}\otimes F'),$$ for any
$x\in W+\gamma$ and $j\in\Z$, and $SS_k(F')\cap(U\times \text{Int}
(\gamma^{\circ a}))=\emptyset$. Hence we may replace $F$ by $F'$ to
prove condition (\ref{E:14}).

Arguing by induction on $k$, we may assume that (\ref{E:14}) holds for $k-1$ and hence $F\in D^{\geq
k}(\mathbf{k}_X)$. Hence, given $x\in W+\gamma$,
 $$H^k(X; \mathbf{k}_{(x+\gamma)\cap H_{\varepsilon}}\otimes
F)\simeq \Gamma(X;  \mathbf{k}_{(x+\gamma)\cap
H_{\varepsilon}}\otimes H^k(F)).$$

Given $s\in \Gamma(X;  \mathbf{k}_{(x+\gamma)\cap
H_{\varepsilon}}\otimes H^k(F))$ we can extend $s$ to a section
$$\widetilde{s}\in \Gamma(\Omega;\mathbf{k}_
{H_{\varepsilon}}\otimes H^k(F))\subset \Gamma(\Omega; H^k(F)),$$
where $\Omega$ is a $\gamma$-open neighborhood of $x+\gamma$ such
that $\Omega\cap H_{\varepsilon}\subset\subset U$.

Set $S=$ supp$(\widetilde{s})\subset \Omega\cap H_{\varepsilon}$.
Since $H^k_{\{\varphi\geq 0\}}(F)\simeq \Gamma_{\{\varphi\geq
0\}}(H^k(F))$, for any real analytic function $\varphi$ defined on
$\R^n$, we get $S=\emptyset$ from the following Lemma, and hence
$H^k(X;  \mathbf{k}_{(x+\gamma)\cap H_{\varepsilon}}\otimes F)=0.$
\end{proof}

\begin{lemma}[~\cite{KMS1}]
Let $\gamma$ be a proper closed convex cone in $\R^n$. Let
$\Omega$ be a $\gamma$-open subset of $\R^n$ and let $S$ be a
closed subset of $\Omega$ such that $S\subset\subset\R^n$. Assume
the following condition: for any $x\in \R^n$ and any real analytic
function $\varphi$ defined on $\R^n$, the three conditions $S\cap
\{x;\varphi(x)< 0\}=\emptyset$, $\varphi(x)=0$ and $d\varphi(x)\in
Int(\gamma^{\circ a})$ imply $x\notin S$. Then $S$ is an empty
set.
\end{lemma}

\begin{corollary}\label{C:1}
Let $E$ be a real finite dimensional vector space, $X$ an open
subset of $E$ and let $F\in D^b(\mathbf{k}_X)$. Let $U$ be an open
subset of $X$ and $\gamma$ be a closed convex proper cone in $E$
with $0\in \gamma$. Assume that $$SS_k(F)\cap(U\times \text{Int}
(\gamma^{\circ a}))=\emptyset.$$ Then, for each $(x_0,\xi_0)\in
U\times \text{Int} (\gamma^{\circ a})$ there exists an open
neighborhood $V$ of $x_0$ in $U$ such that $$R\phi_{\gamma
*}(R\Gamma_{\Omega_1\backslash\Omega_0}( F))\in D^{>
k}(\mathbf{k}_{X_\gamma})$$ for every $\gamma$-open subsets $\Omega_1$ and
$\Omega_0$ with $\Omega_0\subset \Omega_1$,
$\Omega_1\backslash\Omega_0\subset \subset V$ and $x_0\in
\text{Int}(\Omega_1\backslash\Omega_0)$.
\end{corollary}

\begin{proof}[\textbf{Proof}]
We may assume assume $X=\R^n$. Let us consider an open
neighborhood $V'$ of $x_0$ such that $\overline{V'}$ is a compact
subset of $U$.

For each $x\in \overline{V'}$ there exists $\varepsilon(x)>0$ and
an open neighborhood $V(x)$ of $x$ such that $(V(x)+\gamma)\cap
H_{\varepsilon(x)} \subset\subset U$. Since $\overline{V'}$ is
compact, we can find a finite covering of $\overline{V'}$ by open subsets, let us say,
$\{V(x_1),...,V(x_l)\}$, with $l\in\N$. Let us choose
$\varepsilon=min\{\varepsilon(x_i); i=1,...,l \}$. Then,
$(V'+\gamma)\cap H_{\varepsilon} \subset\subset U$, and by Lemma \ref{L:13},
\begin{center}
$H^j(X;\mathbf{k}_{(x+\gamma)\cap H_\varepsilon}\otimes F)=0$, for
each $x\in V'$ and $j\leq k$.
\end{center}

Since $$H^j(X;\mathbf{k}_{(x+\gamma)\cap H_\varepsilon}\otimes F)
\simeq H^j(\phi_\gamma^{-1}R\phi_{\gamma *}
F_{H_\varepsilon})_x,$$ for all $x\in V'$ and $j\leq k$, one gets
$\phi_\gamma^{-1}R\phi_{\gamma *}F_{H_\varepsilon}\in D^{>
k}(\mathbf{k}_{V'_\gamma})$.

Let $V=V'\cap H_\varepsilon$ and consider $\gamma$-open subsets
$\Omega_1$ and $\Omega_0$ such that $\Omega_0\subset\Omega_1$,
$x_0\in \text{Int}(\Omega_1\backslash\Omega_0)$ and
$\Omega_1\backslash\Omega_0\subset \subset V$.

 From $\phi_\gamma^{-1}R\phi_{\gamma *}F_{H_\varepsilon}\in D^{>
k}(\mathbf{k}_{V_\gamma})$ one gets $$R\phi_{\gamma
*}(R\Gamma_{\Omega_1\backslash\Omega_0}(\phi_\gamma^{-1}R\phi_{\gamma
*}F_{H_\varepsilon}))\in D^{> k}(\mathbf{k}_{X_\gamma}).$$ Since
$$R\phi_{\gamma *}(R\Gamma_{\Omega_1\backslash\Omega_0}(
F_{H_\varepsilon}))\simeq R\phi_{\gamma
*}(R\Gamma_{\Omega_1\backslash\Omega_0}(\phi_\gamma^{-1}R\phi_{\gamma
*}F_{H_\varepsilon})),$$ we obtain $R\phi_{\gamma
*}(R\Gamma_{\Omega_1\backslash\Omega_0}( F_{H_\varepsilon}))\in
D^{> k}(\mathbf{k}_{X_\gamma})$.

Finally, we conclude that $R\phi_{\gamma
*}(R\Gamma_{\Omega_1\backslash\Omega_0}( F))\in D^{>
k}(\mathbf{k}_{X_\gamma}),$ since
$R\Gamma_{\Omega_1\backslash\Omega_0}( F_{H_\varepsilon})$ is
isomorphic to $R\Gamma_{\Omega_1\backslash\Omega_0}(F)$.

\end{proof}


\section{Complements on functorial properties of the truncated microsupport}

\hspace*{\parindent} In order to prove the main results we need further
functorial properties of the truncated microsupport similar to
those of the microsupport itself but requiring adapted proofs.

\begin{lemma}\label{P:14i}
Let $X$ be a finite dimensional real vector space, $\gamma$ be a
closed convex proper cone of $X$ with $0\in\gamma$, and $\Omega$ a
$\gamma^a$-open subset of $X$ such that, for any compact $K$ of
$X$, $\Omega \cap (K+\gamma)$ is relatively compact. Let $F\in
D^b(\mathbf{k}_X)$ and assume ${R\phi_{\gamma}} _* F\in
D^{>k}(\mathbf{k}_{X_\gamma})$. Then we have $${R\phi_{\gamma}} _*
F_{\Omega}\in D^{>k}(\mathbf{k}_{X_\gamma}).$$
\end{lemma}

\begin{proof}[\textbf{Proof}]
The proof is contained in the proof of Lemma 5.4.3 (i) of
\cite{K-S1}.
\end{proof}

\begin{proposition}\label{P:14}
Let $X$ be a manifold, $F\in D^b(\mathbf{k}_X)$ and
$\Omega$ be an open subset of $X$.

(i) Assume $SS_k(F)\cap N^*(\Omega)^a\subset T_X^*X$. Then $$SS_k
(R\Gamma _{\Omega} (F))\subset N^*(\Omega) + SS_k(F).$$

(ii) Assume $SS_k(F)\cap N^*(\Omega)\subset T_X^*X$. Then $$SS_k
(F_{\Omega})\subset N^*(\Omega)^a  + SS_k(F).$$
\end{proposition}

\begin{proof}[\textbf{Proof}]
The proof is an adaptation of the proof of Proposition 5.4.8 (i)
and (ii) of~\cite{K-S1}, using Corollary \ref{C:1} and Remark
\ref{R:3} instead of Propositions 5.2.1 and 5.1.1 of~\cite{K-S1},
respectively.
\end{proof}

\begin{proposition}\label{P:15}
Let $\Omega$ be an open subset of $X$ and let $j$ be the embedding
$\Omega\hookrightarrow X$. Let $F\in D^b(\mathbf{k}_\Omega)$.
Then:

(i) $SS_k(Rj_* F)\subset SS_k(F)\widehat{+} N^*(\Omega).$

(ii) $SS_k(Rj_{!} F)\subset SS_k(F)\widehat{+} N^*(\Omega)^a.$
\end{proposition}

\begin{proof}[\textbf{Proof}]
The proof is the stepwise adaptation of the proof of Proposition 6.3.1
of~\cite{K-S1}, using Propositions \ref{P:14}, \ref{PP:1} and
Corollary \ref{C:1} instead of Propositions 5.4.8, 5.1.1 and 5.2.1
of~\cite{K-S1}, respectively.
\end{proof}

\section{Proofs of the main results}

\subsection{Proofs of Theorems \ref{T:26},
 \ref{P:6} and Corollaries}

\textbf{\emph{Proof of Theorem \ref{P:6}}} Let us first consider
the case of the embedding of a closed submanifold of $X$:

\begin{proposition}\label{L:30}
Let $M$ be a closed submanifold of $X$ and $F\in
D^b(\mathbf{k}_X)$. Then $$SS_k(F|_M)\subset
j_dj_\pi^{-1}(SS_k(F)\widehat{+}T_M^*X),$$ where $j$ is the
embedding of $M$ in $X$.
\end{proposition}

\begin{proof}[\textbf{Proof}]
Let $d$ denote the codimension of $M$.
Let $(x_1,...,x_n)$ be a system of local coordinates on $X$ such
that $M=\{(x_1,...,x_n); x_1=...=x_d=0\}$ and let $(x;\xi)$ denote the
associated coordinates on $T^*X$. Set $x'=(x_1,...,x_d)$, $x''=(x_{d+1}, ..., x_{n})$.

Let $(x''_0;\xi''_0)\in T^*M$ such that $(x''_0;\xi''_0)\notin
j_dj_\pi^{-1}({SS_k(F)\widehat{+}T_M^*X})$. We shall prove that
$(x''_0;\xi''_0)\notin SS_k(F|_M)$.

By the assumption, $(0,x''_0;\xi',\xi''_0)\notin
{SS_k(F)\widehat{+}T_M^*X}$ for any $\xi'\in \R^d$. In particular,
$(0,x''_0;0,\xi''_0)\notin {SS_k(F)\widehat{+}T_M^*X}$. We may
assume that $(0,x''_0)\in \pi(SS_k (F))\cap M$ and by (\ref{tau}),
that $\xi''_0\neq 0$.

Setting $x_0=(0,x''_0)$, $\xi_0=(0,\xi''_0)$ and $p=(x_0,\xi_0)$,
there exists a closed convex proper cone $\gamma$ such that Int$(\gamma)\neq \emptyset$ and
\begin{equation}\label{E12}
\begin{cases} \xi_0\in\text{Int}(\gamma^{\circ a}),\\ (\{x_0\}\times
\gamma^{\circ a})\cap (SS_k(F)\widehat{+}T_M^*X)\subset T^*_{X} X.
\end{cases}
\end{equation}
Therefore we may find  a
neighborhood $V$ of $x_0$ such that
\begin{equation}\label{E13}
(V\times \gamma^{\circ a})\cap( SS_k(F)\widehat{+}T_M^*X) \subset
T^*_{X} X.
\end{equation}
In particular,  $(\{x_0\}\times\gamma^{\circ a})\cap
SS_k(F)\subset \{(x_0;0)\}$. Therefore, $$(\{x_0\}\times
\text{Int}(\gamma^{\circ a}))\cap SS_k(F)=\emptyset,$$ and we may
choose $V$ such that
\begin{equation}\label{E113}
(V\times \text{Int}(\gamma^{\circ a}))\cap SS_k(F)=\emptyset.
\end{equation}
and
\begin{equation}\label{E114}
(V\times \gamma^{\circ a})\cap T^*_M  X\subset T^*_X X.
\end{equation}

After changing the local coordinates on $X$ if necessary, we may
also assume: $$\begin{cases} \xi_0=(0,...,0,1),\\ \gamma^{\circ
a}=\{(\xi',\xi''); \xi_n\geq \delta
\vert(\xi',\xi_{d+1},...\xi_{n-1}) \vert \}, \end{cases}$$ for
some $\delta> 0$. Hence, $$\gamma=\{(x',x''); x_n\leq
-\frac{1}{\delta} \vert(x',x_{d+1},...x_{n-1}) \vert\},$$ and, for
any $x\in M, $ $(x+\gamma)\cap M=x+(\gamma\cap M)$. For
$\varepsilon>0$ let us denote by $H_\varepsilon$ the open
half-space $ H_\varepsilon=\{x\in X; \text{Re}\langle x-x_0,
\xi_0\rangle >-\varepsilon\}$. Let us choose $\varepsilon> 0$ and
an open neighborhood $W\subset V$ of $x_0$ such that
$(W+\gamma)\cap H_\varepsilon\subset\subset V$. Set
$\gamma'=\gamma\cap M$, $V'=V\cap M$,  $W'=W\cap M$ and
$H'_\varepsilon=\{x''\in M; \langle x''-x''_0,\xi''_0\rangle>
-\varepsilon\}$. Since $\gamma'$ is a closed convex proper cone in
$M$ such that $\xi'_0\in\text{Int}(\gamma'^{\circ a})$ and $W'$ is
an open neighborhood of $x'_0$ in $M$,  by Proposition \ref{PP:1}
its enough to prove that there exists $\varepsilon' >0 $ such that
$(W'+\gamma')\cap \overline{H'_{\varepsilon'}}\subset M$ and $H^j
(M; \mathbf{k}_{(x+\gamma')\cap H'_{\varepsilon'}}\otimes
F|_M)=0,$ for all $j\leq k$ and $x\in W'$.

This will be a consequence of Lemma \ref{L:131} with
$\Lambda=SS_k(F)$. We shall use the notation $\gamma_{\lambda},
V_{\lambda}$ introduced in Section 3. Let $C$ be given by Lemma
\ref{L:131} and let us choose sequences $(\lambda_l)_{l
\in\mathbb{N}}$, $(\epsilon_l)_{l\in\mathbb{N}}$ of real positive
numbers, satisfying $0<\epsilon_l,\lambda_l<C$, such that
$(\lambda_l)_{l\in\mathbb{N}}$ converges to $0$ and
$(\epsilon_l)_{l\in\mathbb{N}}$ converges to $C$.

Remark that $\gamma_{\lambda_{l}}^{\circ a}\supset \gamma^{\circ
a}$ and that $$V_{\lambda_{l} \epsilon_{l}}\cap M=V\cap M=V',$$
$$W_{\lambda_{l} \epsilon_{l}}\cap M=W\cap M=W',$$
$$W_{\lambda_{l} \epsilon_{l}} +\gamma_{\lambda_{l}
\epsilon_{l}}\cap H_{\varepsilon}\subset V_{\lambda_{l}
{}\epsilon_{l}}, $$ $$(W_{\lambda_{l} \epsilon_{l}}+\gamma
_{\lambda_{l}})\cap H_{\epsilon_l \varepsilon}\subset
V_{\lambda_{l} {}\epsilon_{l}}.$$ Let $x''\in M\cap W$ be given,
choose a sequence $x''_{l}$ in $W$ converging to $x''$ such that
$x''\in \text{Int} (x''_{l}+ \gamma')$ and note $H'=H\cap M$.
Then, for any $j\geq k$, we have $$H^j (M;
\mathbf{k}_{(x''+\gamma')\cap H'_{C \varepsilon}}\otimes F\vert_
M) \simeq \lim_{\substack{\longrightarrow\\ l}} H^j (X;
\mathbf{k}_{(x''_{l}+\gamma_{\lambda_{l}})\cap H_{\epsilon_l
\varepsilon}}\otimes F)=0,$$ thanks to Lemma \ref{L:13}. Hence
$(x''_{0};\xi''_{0})\notin SS_{k}(F\vert _{M})$.
\end{proof}

\begin{proof}[\textbf{End of the proof of Theorem \ref{P:6}}]Let us
decompose $f$ by the graph map
\begin{center}
$Y\xrightarrow[g]{}Y\times X\xrightarrow[h]{}X$, $f=h\circ g$
\end{center}
where $g(y)=(y,f(y))$ and $h$ is the second projection on $Y\times
X$.

Identifying $Y$ with the graph of $f$, we may assume that $Y$ is a
closed subvariety of $Y\times X$, and we get by Proposition \ref{L:30}
 and Proposition \ref{P:3},
$$SS_k(f^{-1}F)=SS_k(g^{-1}(h^{-1}F))\subset $$ $$\subset
g_dg_\pi^{-1}(h_dh_\pi^{-1}(SS_{k}(F))\widehat{+}T_Y^*(Y\times
X)).$$

We shall prove that
$$g_dg_\pi^{-1}(h_dh_\pi^{-1}(SS_{k}(F))\widehat{+}T_Y^*(Y\times
X))= f^\#(SS_k(F)).$$

Let $(y)$ be a system of local coordinates  on $Y$, $(x)$ a system
of local coordinates on  on $X$ and let $(y;\xi)$, $(x;\eta)$ be
the associated coordinates on $T^*Y$ and $T^*X$, respectively.

Let $(y_0;\xi_0)\in
g_dg_\pi^{-1}(h_dh_\pi^{-1}(SS_{k}(F))\widehat{+}T_Y^*(Y\times
X)),$ then there exists $\xi, \eta$ such that
$(y_0,f(y_0);\xi,\eta)\in
h_dh_\pi^{-1}(SS_{k}(F))\widehat{+}T_Y^*(Y\times X)$ and
$\xi_0=\xi+{^t}f'(y_0)\cdot\eta$. Hence we may find sequences
$\{(y_n,x_n;\xi_n,\eta_n)\}_n$ in $h_dh_\pi^{-1}(SS_{k}(F))$ and
$\{(y'_n,f(y'_n);\xi'_n,\eta'_n)\}_n$ in $T_Y^*(Y\times X)$ such
that $$\begin{cases} (y_n,x_n), (y'_n,f(y'_n)) \xrightarrow[n]{}
(y_0, f(y_0)),\\ (\xi_n,\eta_n)+(\xi'_n,\eta'_n) \xrightarrow[n]{}
(\xi,\eta),\\
|(y_n,x_n)-(y'_n,f(y'_n))||(\xi_n,\eta_n)|\xrightarrow[n]{} 0.
\end{cases}$$

One has $(x_n;\eta_n)\in SS_k(F),$ $\xi_n=0$ and
$\xi'_n+{^t}f'(y'_n)\cdot\eta'_n=0$; hence we obtain
${^t}f'(y'_n)\cdot(\eta_n+\eta'_n)\xrightarrow[n]{}{^t}f'(y_0)\cdot\eta=\xi_0-\xi$
and then ${^t}f'(y'_n)\cdot\eta_n\xrightarrow[n]{}\xi_0.$

Therefore we have sequences $\{(x_n;\eta_n)\}_n\in SS_k(F)$ and
$\{y'_n\}_n$ in $Y$ such that $$\begin{cases} y_n\xrightarrow[n]{}
y_0, x_n\xrightarrow[n]{} f(y_0),\\
{^t}f'(y'_n)\cdot\eta_n\xrightarrow[n]{}\xi_0\\
|x_n-f(y'_n)||\eta_n|\xrightarrow[n]{} 0.
\end{cases}$$

This gives $(y_0;\xi_0)\in f^\#(SS_k(F))$ and also the converse
thanks to Proposition \ref{P:132}.
\end{proof}

\begin{corollary}
Let $M$ be a closed submanifold of $X$ and $F\in
D^b(\mathbf{k}_X)$. Then $$SS_k(F_M)\subset
SS_k(F)\widehat{+}T_M^*X.$$
\end{corollary}

\begin{proof}[\textbf{Proof}]
Let $j:M\hookrightarrow X$ denote the embedding of $M$ on $X$.
Then $F_M\simeq j_*(F|_M)$ and by Proposition \ref{P:2} and
Proposition \ref{L:30} $$SS_k(F_M)=j_\pi
j_d^{-1}(SS_k(F|_M))\subset$$ $$\subset j_\pi j_d^{-1}j_d
j_\pi^{-1}(SS_k(F)\widehat{+}T_M^*X)\subset
SS_k(F)\widehat{+}T_M^*X.$$
\end{proof}

\begin{proposition}\label{P:16}
Let $M$ be a closed submanifold of $X$, $U=X\setminus M$, $j$ the
embedding $U\hookrightarrow X$, $\iota$ the embedding of $M$ in
$X$ and let $F\in D^b(\mathbf{k}_U)$. Then:

(i) $SS_k(Rj_* F)\cap \pi^{-1} (M)\subset SS_k(F)\widehat{+}
T_M^*X,$

(ii) $SS_k(Rj_! F)\cap \pi^{-1}(M)\subset SS_k(F)\widehat{+}
T_M^*X,$

(iii) $SS_k((Rj_*F)|_M)\subset
\iota_d\iota_\pi^{-1}(SS_k(F)\widehat{+}T_M^*X).$
\end{proposition}

\begin{proof}[\textbf{Proof}]
The proof of the two first conditions is analogous to the proof of
the two first conditions of Proposition 6.3.2 of ~\cite{K-S1},
replacing Proposition 5.4.4 and Theorem 6.3.1 of ~\cite{K-S1}, by
Propositions \ref{P:2} and \ref{P:15}, respectively.

Let us now prove the third inequality. By Proposition \ref{L:30} and $(i)$,

$$ SS_k((Rj_*F)|_M)\subset
\iota_d\iota_\pi^{-1}(SS_k(Rj_*F)\widehat{+}T_M^*X)\subset$$
$$\subset \iota_d\iota_\pi^{-1}((SS_k(F)\widehat{+}
T_M^*X)\widehat{+}T_M^*X).$$

By Lemma \ref{L:133} $$(SS_k(F)\widehat{+}
T_M^*X)\widehat{+}T_M^*X= SS_k(F)\widehat{+}T_M^*X.$$ Hence
$$SS_k((Rj_*F)|_M)\subset
\iota_d\iota_\pi^{-1}(SS_k(F)\widehat{+}T_M^*X).$$
\end{proof}

Note that, with Lemma \ref{L:32} and Proposition \ref{P:16} in hand, we
obtain the  analogue of Proposition 6.3.2 of \cite{K-S1}.

\begin{corollary}\label{C:2}
Let $M$ be a closed submanifold of $X$ and $F\in
D^b(\mathbf{k}_X)$. Then $$SS_k(R\Gamma_{M}(F))\subset
SS_k(F)\widehat{+} T_M^*X.$$
\end{corollary}

\begin{proof}[\textbf{Proof}]
This is a consequence of Proposition \ref{P:16}, together with the
distinguished triangle $$R\Gamma_M(F)\rightarrow F\rightarrow
R\Gamma_{X\backslash M}(F) \xrightarrow{+1}.$$
\end{proof}

\begin{corollary}\label{C:22}
Let $M$ be a closed submanifold of $X$ and $F\in
D^b(\mathbf{k}_X)$. Assume that $$SS_k(F)\cap T_M^*X\subset
T_X^*X.$$ Then we have a natural isomorphism $$\tau^{\leq
k}(F_M\otimes \omega_{M|X})\simeq \tau^{\leq k}(R\Gamma_M(F)).$$In
particular $$SS_k(F_M\otimes \omega_{M|X})=SS_k(R\Gamma_M(F)).$$
\end{corollary}

\begin{proof}
Let $\dot{\pi}$ be the restriction of $\pi$ to the cotangent
bundle deprived of the zero section. We have a distinguished
triangle $$\tau^{\leq k}(F_M\otimes \omega_{M|X})\to \tau^{\leq k}(R\Gamma_M(F))\to
\tau^{\leq k}(\dot{\pi}(\mu_M (F)|_{\dot{T}^*_M X}))\underset{+1}{\to}.$$ Since
by Theorem 5.1 of \cite {KMS1} $supp(\tau^{\leq
k}(\mu_M(F)))\subset SS_k(F)\cap T_M^*X$,  the
third term vanishes hence the result.
\end{proof}

\begin{corollary}\label{C:10}
Let $M$ be a closed submanifold of $X$ and $F\in
D^b(\mathbf{k}_X)$. Let $j$ denote the embedding of $M$ in
$X$.Then $$SS_k(j^{!}F)\subset j_d j_{\pi}
^{-1}(SS_k(F)\widehat{+} T_M^*X).$$
\end{corollary}

\begin{proof}[\textbf{Proof}]
This is a consequence of Corollary 6.4 and Proposition \ref{L:30}
together with Lemma \ref{L:133}.
\end{proof}

\begin{proof}[\textbf{Proof of Theorem \ref{T:26}}]
The proof is the adaptation step by step of the proof of Theorem 6.4.1 of
~\cite{K-S1}, applying Proposition \ref{P:3},  and Proposition
\ref{P:16} instead of Proposition 5.4.5 and Proposition 6.3.2,
respectively of~\cite{K-S1}.
\end{proof}

Let now $Y$ be a complex closed smooth hypersurface of $X$ defined
as the zero locus of a holomorphic function $f$. Let $\psi_Y$
denote the functor of nearby cycles associated to $Y$. Then $Y$
may be regarded as a submanifold of $T_Y X$ by a canonical section
$s$ such that $\psi_Y (F)\simeq s^{-1}\nu_Y (F)$.  Once more we identify
$T_{T_Y^*X}T^*X$, $T^*(T^*_Y X)$ and $T^*(T_Y X)$ (cf. Proposition
5.5.1 of \cite{K-S1}).

Recall that, in a system of linear coordinates $x=(x_1,...,x_n)$
on $X$ such that $Y$ is defined by $x_1=0$, $s\colon Y\to T_Y X$
is the section $s(x_2,...,x_n)=(x_2,...,x_n; 1)$.  With the local
coordinates described above, and $A$ being a conic closed subset
of $T^*(T^*_Y X)$, we have: $$s_d s_\pi^{-1}(A)= \{(x_2,...,x_n;
\xi_2,...,\xi_n); \exists \xi_1, ( x_2,...,x_n,
1;\xi_2,...,\xi_n,\xi_1)\in A\}.$$

Corollary 1.4 is an  immediate consequence of Proposition
\ref{L:30}.

The following estimate for the tensor product can be seen as a
generalization of Proposition \ref{L:30}:

\begin{proposition}\label{P:28}
Let $F$ and $G$ belong to $D^b(\mathbf{k}_X)$. Then:
$$SS_k(F\otimes^{\mathbb{L}} G)\subset
\underset{i+j=k}{\bigcup}(SS_i(F)\widehat{+}SS_j(G)).$$
\end{proposition}

\begin{proof}[\textbf{Proof}]
Let $\delta_X:X\rightarrow X\times X$ be the diagonal embedding.

Since $F\otimes^{\mathbb{L}} G\simeq \delta_X^{-1}(F\boxtimes^{\mathbb{L}} G)$, the result
follows from Proposition \ref{L:30} and Proposition \ref{P:21}.
\end{proof}

\subsection{Application to $\shd$-modules}\label{S:2}

\hspace*{\parindent} Let $X$ be a complex finite dimensional manifold. One of the
important problems in the theory of $\cal{D}$-modules is the
relation between the caracteristic variety of a system $\cal{M}$
and that of its induced system ${\cal{M}}_{Y}$ along a closed
submanifold $Y$, which was completely solved in the non
charateristic case by M. Kashiwara as well as in a more general
situation treated in \cite{LS}, which includes the case where
$\cal{M}$ is regular along $Y$ in the sense of \cite{KO}.
Similarly, in the case of a smooth complex hypersurface,  it is
interesting to relate $\text{Char} (\cal{M})$ and
$\text{Char}(\psi_Y (\cal{M}))$, where $\psi_Y$ denotes the functor
of nearby cycles.

Let $d$ be the codimension of $Y$, denote by $j$ the embedding
$Y\to X$ and by $\pi'$ the projection $T^* Y \to Y$. Given an
homogeneous involutive subvariety $V$ of $T^*X$ of codimension
$\geq d$, we shall say that $Y$ is orthogonal to $V$ if there
exists a smooth involutive submanifold $V^*$ containing $V$ such
that $Y$ and $V^*$ are orthogonal. More precisely, there exist a
set $\{f_1,..., f_d\}$ of homogeneous functions of degree zero
vanishing on $\pi^{-1}(Y)$, such that the differential $df_i$ are
linearly independent on $\pi^{-1}(Y)$, and a set  $\{g_1,...,
g_p\}, p\geq d,$ of homogeneous functions of degree one linearly
independent on $V^*$ such that the matrix of the Poisson brackets
$[\{f_i,g_j\}]|_{V^*}$ has everywhere rank $d$.

As before, $F$ will denote the complex
$R\mathcal{H}{om}_{\shd_X}(\mathcal{M},\sho_X)$. Let
$SS(F)=\bigsqcup_\alpha V_\alpha$ be the decomposition of $SS(F)$
in its irreducible involutive components in a neighborhood of
$p\in T^*X$. Let us denote by $Y_{\alpha}$ the variety $\pi
(V_{\alpha})$.

Recall that in Theorem 6.7 of \cite{KMS1} it is proved that, for
any $k$, $SS_k(F)=(\underset{codim Y_{\alpha}< k}{\bigcup}
V_\alpha)\cup(\underset {codimY_{\alpha}= k}{\bigcup}
T^*_{Y_{\alpha}}X)$.

\begin{proof}[\textbf{Proof of Theorem \ref{P:146}}]
The first assertion is an immediate consequence of Theorem
\ref{T:26} and the second follows from the regularity of $\shm$.
\end{proof}

\begin{proof}[\textbf{Proof of Corollary \ref{P:45}}]
It is a consequence of Corollary \ref{P:44} and the regularity of $\shm$.
\end{proof}
\begin{proof}[\textbf{Proof of Theorem \ref{P:46}}] Since $\shm$ is regular along
$Y$, one has the isomorphism
$$R\mathcal{H}{om}_{\shd_Y}(\mathcal{M}_Y,\sho_Y)\simeq
R\mathcal{H}{om}_{\shd_X}(\mathcal{M},\sho_X)|_Y,$$ and the first
part is an immediate consequence of Proposition \ref{L:30}.
Let us now prove the second assertion. It will be a consequence of
the Lemma below:

\begin{lemma}\label{L:28}
Assume that the homogeneous involutive variety $V$ is irreducible
and that $Y$ is orthogonal to $V$.

Then :

(i) $V'=j_d(j_\pi^{-1}(V))$ is an irreducible homogeneous
involutive subvariety of $T^*Y$.

(ii) The
codimension of $\pi'(V')$ is equal to the codimension of
$\pi(V)$.

(iii) When $V$ is the characteristic variety of a
coherent $\shd_X$-module, the orthogonality of $Y$ implies that
$Y$ is non characteristic for $\shm$.
\end{lemma}

\begin{proof}[\textbf{Proof}]
Let $V^*$ be a smooth involutive manifold containing $V$ such that
$Y$ is orthogonal to $V^*$. Since the assertions can be checked
locally, by a standard reasoning we may consider a system  $(x;
\xi)$  of local symplectic coordinates on $T^*X$ in a neighborhood
of $p\in V\cap \pi^{-1}(Y)=j_\pi^{-1}(V)$, such that $ Y$ is the
sumanifold $\{(x)=(x_1,...x_n); x_1=...=x_d=0\}$ and $V^*$ is
defined in $T^*X$ by the equations $\xi_1=...=\xi_d= g_{d+1}(x'';
\xi'')= ...=g_p(x''; \xi'')=0$, where we set  $(x')= (x_,...,
x_d)$ (resp. $(\xi')=(\xi_1,...\xi_d)$), $(x'')=
(x_{d+1},...,x_{n})$ (resp. $(\xi'')= (\xi_{d+1},...,\xi_{n})$).
Therefore, the irreducible ideal of definition $I(V)$ is generated
by a set of functions $$\{\xi_1,..., \xi_d, g_{d+1}(x'';
\xi''),...,g_p(x''; \xi''), h_{p+1}(x";\xi"),...,
h_{p+l}(x";\xi")\},$$  for some $l\geq 0$. Hence $I(V')$ is
generated in $\sho_{T^*Y}$ by the set of functions $$\{
g_{d+1}(x''; \xi''), ..., g_p(x''; \xi''), h_{p+1}(x";\xi"),...,
h_{p+l}(x";\xi")\},$$ which entails $(i)$, $(ii)$ and $(iii)$.
\end{proof}

Since $Y$ is non characteristic, we have
$$SS(F|_Y)=\text{Char}(\shm_Y)=j_d j_{\pi}^{-1}(SS(F)).$$

On the other side, since $SS_k(F)\cap  T^*_Y X\subset T^*_X X$, we
get from the first assertion that $SS_k(F|_Y)\subset j_d
j_{\pi}^{-1}(SS_k(F))$, for any k. Moreover, setting
$V'_{\alpha}=j_d j_\pi^{-1} (V_{\alpha})$, by the preceding Lemma,
for any $\alpha$ such that $codim Y_{\alpha}\leq k$, $V'_\alpha$
is an irreducible component of $SS(F|_Y)$. Therefore by Theorem
6.7 of \cite {KMS1}, for any  $i\leq k,$  $$ SS_i(F|_Y)\supset j_d
j_{\pi}^{-1}(SS_i(F)).$$
\end{proof}

\begin{example}
Let $X=\mathbb{C}^n$, with $n\geq 3$, endowed with the coordinates
$(x_1,...,x_n)$. Let $Y$ be the hypersurface $\{x_n=0\}$ and
$\Omega=\{x\in X; \text{Re}( x_1-x_{n-1})< 0\}$. Let $\mathcal{J}$
be a coherent left ideal of $\shd_X$ and set
$\shm=\shd_X/\mathcal{J}$. Assume that there exist  in
$\mathcal{J}$ an operator $P$ in the Weierstrass form with respect
to the derivation $D_{x_n}$ and an operator $Q$ such that the
principal symbol of $Q$, $\sigma (Q)$, is of the form
$$\sigma(Q)=x_1 q(x_1,..., x_{n-1};\xi_1,...,\xi_{n-1}),$$ and $q$
does not vanish on $T^*_{\delta\Omega}X$. Then,$ T^*_{\delta(\omega)}X\cap SS_1(\shm)\subset \{0\}$ and, setting
$\Omega'=\Omega\cap Y$, $\Omega'$ has smooth boundary
$\delta\Omega'$. By Theorem 1.3, $T^*_{\delta\Omega'} Y\cap
SS_1(\shm_Y)\subset\{0\}.$ Therefore
$$\mathcal{H}\text{om}_{\shd_Y}(\shm_Y,
\mathcal{H}^1_{\{\text{Re}(x_1-x_{n-1})\geq 0\}}
(\sho_Y))|_{\delta\Omega'}=0.$$
\end{example}

\begin{proof}[\textbf{Proof of Corollary 1.5}]
As proved in \cite{K5}, we have the isomorphism
$$R\mathcal{H}{om}_{\shd_Y}(\psi_Y (\shm),\sho_Y)\simeq
\psi_Y(R\mathcal{H}{om}_{\shd_X}(\mathcal{M},\sho_X)).$$ It is
then enough to use Proposition \ref{L:30}.
\end{proof}

Let $M$ be a real analytic manifold of dimension $n$, $X$ a
complex analytic manifold complexifying $M$ and $\mathcal{M}$ a
coherent $\shd_X$-module.

Let $\sha_M$ denote the sheaf of real analytic functions on $M$.
Remark that $\sha_M=\sho_X|_{M}$. Let  $\shb_M$
denote the sheaf of Sato's hyperfunctions on $M$. Recall that
$$\shb_M\simeq R\Gamma_M(\sho_X)\otimes {or}_{M/X}.$$

\begin{proof}[\textbf{Proof of Proposition \ref{P:7}}]
One has $$ R\mathcal{H}{om}_{\shd_X}(\mathcal{M},\sha_M)\simeq
F|_M.$$

Therefore, by Proposition \ref{L:30}
$$SS_{k}(R\mathcal{H}{om}_{\shd_X}(\mathcal{M},\sha_M))\subset
j_dj_\pi^{-1}(SS_k(F)\widehat{+}T_M^*X)$$
\end{proof}
Let us remark that a variant of the preceding result was obtained in \cite{T}
for $k=1$ using directly the properties of holomorphic functions.

\begin{proof}[\textbf{Proof of
Corollary \ref{P:118}}]
The first part is an immediate consequence of Corollary \ref{C:22}. The second follows from Proposition \ref{P:7} and Theorem 6.7 of \cite{KMS1}.
\end{proof}

\begin{example}\label{eg}
Let $M=\mathbb{R}^n$, with $n\geq 2$, endowed with the coordinates
$x=(x_1,...,x_n)$. Let $\Omega=\{x\in M; \phi(x)< 0\}$ for some
real $C^1$-function. Let $X=\mathbb{C}^n$ and $\shm$ be a coherent
$\shd_X$-module defined by $\shm=\shd_X/\shd_X P$ where $P$ is a
differential operator. Assume that the principal symbol  $\sigma
(P)$ is of the form $$\sigma(P)=a(x) q(x;\xi),$$ where $a(x)$ is
an holomorphic function and $q$ does not vanish on $T^*M$, more
precisely, $q$ is the principal symbol of an elliptic operator.
Recall that $ SS_1(F)\subset \overline{\{(x;\xi); a(x)=0,
\xi\in\mathbb{C}da(x)\}}\cup q^{-1} (0).$This entails that $T^*_M
X\cap SS_1(F)\subset T^*_X X$ hence
$$SS_{1}(R\mathcal{H}{om}_{\shd_X}(\mathcal{M},\sha_M))=SS_{1}(R\mathcal{H}{om}_{\shd_X}(\mathcal{M},\shb_M))\subset
j_d j_{\pi}^{-1}(SS_1(F)).$$

Assume that $d\phi(x)$ is not in $\overline{\mathbb{C}da(x)}$ for
any $x\in\delta \Omega\cap a^{-1}(0)$.    Hence
$T^*_{\delta\Omega} M\cap
SS_1(R\mathcal{H}{om}_{\shd_X}(\mathcal{M},\sha_M))\subset T^*_M
M.$ In other words $$\mathcal{H}\text{om}_{\shd_X}(\shm,
\mathcal{H}^1_{\{\phi(x)\geq 0\}} (\sha_M))|_{\delta\Omega}$$  $$=\mathcal{E}\text{xt}^1_{\shd_X}(\shm, \Gamma_{\{\phi(x)\geq 0\}}(\shb_M))|_{\delta\Omega}=0.$$
\end{example}

\begin{remark}
In general we do not have an interesting estimate for \\
$SS_k(R\mathcal{H}{om}_{\shd_X}(\mathcal{M},\shb_M))$. Let $j$
denote the inclusion $M\hookrightarrow X$. Then $\shb_M\simeq
j^!\sho_X\otimes {or}_{M/X}$ and
$R\mathcal{H}{om}_{\shd_X}(\mathcal{M},\shb_M)\simeq
j^!(R\mathcal{H}{om}_{\shd_X}(\mathcal{M},\sho_X))[n]$. By
Corollary \ref{C:10}, one gets
$$SS_k(R\mathcal{H}{om}_{\shd_X}(\mathcal{M},\shb_M))=
SS_{k+n}(j^!R\mathcal{H}{om}_{\shd_X}(\mathcal{M},\sho_X))\subset$$
$$\subset
SS_{k+n}(R\mathcal{H}{om}_{\shd_X}(\mathcal{M},\sho_X))\widehat{+}T_M^*X.$$
By Theorem 6.7 of ~\cite{KMS1},
$$SS_{k+n}(R\mathcal{H}{om}_{\shd_X}(\mathcal{M},\sho_X))=SS(R\mathcal{H}{om}
_{\shd_X}(\mathcal{M},\sho_X))= \text{Char}(\cal{M}).$$

Hence we get
$$SS_k(R\mathcal{H}{om}_{\shd_X}(\mathcal{M},\shb_M))\subset
\text{Char}(\mathcal{M})\widehat{+}T_M^*X  \ \text{for any k}\geq 0,$$ in
other words, if $M$ is hyperbolic for $\mathcal{M}$ then
$$SS_{k+n}(R\mathcal{H}{om}_{\shd_X}(\mathcal{M},\sho_X))\subset
T_X^*X.$$ But this is well known and is an example that the notion
of truncated microsupport does not work well under Fourier
Transform.
\end{remark}

Let $D^b_{\C-c}(\mathbf{k}_X)$ denote the full subcategory of
$D^b(\mathbf{k}_X)$ consisting of objects with $\C$-constructible
cohomology, that is, the objects $F\in D^b(\mathbf{k}_X)$ for
which there exists a complex analytic stratification $X=\bigcup
X_\alpha$ such that the sheaf $H^j(F)|_{X_\alpha}$ is locally
constant of finite rank, for every $j\in\Z$ and $\alpha$.

A perverse sheaf is an object $F$ of $D^b_{\C-c}(\mathbf{k}_X)$
satisfying the  following two conditions:

\emph{(a) for any complex submanifold $Y$ of $X$ of codimension
$d$, $H^j_Y(F)|_Y$ is zero for $j< d$;}

\emph{(b) for any $j\in\Z$, $H^j(F)$ is supported by a complex
analytic subset of codimension $\geq j$.}

P. Schapira  proved in ~\cite{S2} that, when $F$ is a perverse
object of $D^b_{\C-c}(\mathbf{k}_X)$,
\begin{center}
$H^j(R\Gamma_S(F))_x=0$, for $j\geq 2n$,
\end{center}
for any closed subanalytic subset $S$ of $X$ and any $x\in X$
being non isolated in $S$.

\begin{proposition}\label{P:50}
Let $\mathcal{M}$ be a coherent $\shd_X$-module
Then
$$SS_{n-1}(R\mathcal{H}{om}_{\shd_X}(\mathcal{M},\shb_M))=
SS(R\mathcal{H}{om}_{\shd_X}(\mathcal{M},\shb_M)).$$
\end{proposition}
\begin{proof}[\textbf{Proof}]
Let $\varphi$ be a real analytic function defined on $X$ and
$x_0\in X$ such that $\varphi(x_0)=0$. Then the set $\{x\in X;
\varphi(x)\geq 0\}$ is a closed subanalytic subset of $X$. Assume
that $\cal{M}$ is holonomic. By the Riemann-Hilbert correspondence
(~\cite{K4}),  $F$  is perverse.

Hence, $$H^j(R\Gamma_{\{\varphi\geq
0\}}R\mathcal{H}{om}_{\shd_X}(\mathcal{M},\shb_M))_{x_0}\simeq
H^{j+n}(R\Gamma_{\{\varphi\geq 0\}\cap M}(F))_{x_0}=0,$$ for every
$j\geq n$. By Proposition \ref{PP:1},
$$SS_{n-1}(R\mathcal{H}{om}_{\shd_X}(\mathcal{M},\shb_M))=
SS(R\mathcal{H}{om}_{\shd_X}(\mathcal{M},\shb_M)),$$ under the
assumption that $\mathcal{M}$ is an holonomic $\shd_X$-module.

To treat the general case, we argue as in the proof of
Theorem 2 of ~\cite{S2}. Let us denote by $^*$ the functor
$\mathcal{N}\rightarrow
\mathcal{N}^*=\mathcal{E}{xt}^n_{\shd_X}(\mathcal{N},\shd_X).$ Recall that
Kashiwara  proved in ~\cite{K2} that if $\mathcal{M}$ is coherent, then
$\mathcal{M}^*$ is holonomic,
$\mathcal{M}^{***}\simeq\mathcal{M}^*$ and $\mathcal{M}^{**}$ is a
submodule of $\mathcal{M}$. Defining the coherent $\shd_X$-module
$\mathcal{L}$ by the exact sequence:
$$0\rightarrow\mathcal{M}^{**}\rightarrow\mathcal{M}\rightarrow\mathcal{L}\rightarrow
0,$$ one gets $\mathcal{L}^*=0$ and so $\mathcal{L}$ locally
admits a projective resolution of lenght $n-1$. Therefore,
$$H^j(R\Gamma_{\{\varphi\geq
0\}}R\mathcal{H}{om}_{\shd_X}(\mathcal{M},\shb_M))_{x_0}\simeq
H^j(R\Gamma_{\{\varphi\geq
0\}}R\mathcal{H}{om}_{\shd_X}(\mathcal{M}^{**},\shb_M))_{x_0}=0.$$
for $j\geq n$.

This proves
$$SS_{n-1}(R\mathcal{H}{om}_{\shd_X}(\mathcal{M},\shb_M))
=SS(R\mathcal{H}{om}_{\shd_X}(\mathcal{M},\shb_M)),$$
for every coherent $\shd_X$-module $\mathcal{M}$.
\end{proof}

{\small Ana Rita Martins\\ Centro de {\'A}lgebra da Universidade
de Lisboa, Complexo 2,\\ 2 Avenida Prof. Gama Pinto, 1699 Lisboa
 Portugal\\ arita@pcmat.fc.ul.pt\newline

Teresa Monteiro Fernandes\\ Centro de {\'A}lgebra da Universidade
de Lisboa, Complexo 2,\\ 2 Avenida Prof. Gama Pinto, 1699 Lisboa
 Portugal\\ tmf@ptmat.fc.ul.pt}


\begin{thebibliography}{15}
\bibitem{BS}
Bony J. M. and Schapira P., {\em Existence et prolongement des
solutions holomorphes des \' Equations aux d\' eriv\' ees
partielles}, Inventiones Math. \textbf{17}  95-105, (1972).


\bibitem {E-K-S}
Ebenfelt P., Khavinson D. and Shapiro H., {\em Extending solutions
of holomorphic partial differential equations across real analytic
hypersurfaces}, J.London Math.Soc. (2) \textbf{57}, 411-432,
(1998).

\bibitem{Ho}
H\"ormander L., {\em The analysis of linear partial differential
operators II,} Grundlehren der Math. Wiss.\,\textbf{257}, Springer
Verlag, (1983).


\bibitem{K3}
Kashiwara M., {\em Algebraic Analysis,},  AMS (2001).

\bibitem{K1}
Kashiwara M., {\em On the maximally overdetermined systems of
linear differential equations, I,} Publ. RIMS, Kyoto Univ.,
\textbf{10}, 563-579, (1975).

\bibitem{K2}
Kashiwara M., {\em b-functions and holonomic systems,} Inventiones
Math. \textbf{38}, 33-53, (1976).

\bibitem{K4}
Kashiwara M., {\em The Riemann-Hilbert problem for holonomic
systems,} Publ. RIMS, Kyoto Univ., \textbf{20}, 319-365,(1984).

\bibitem{K5}
Kashiwara M., {\em Vanishing cycle sheaves and holonomic systems
of differential equations,} Lect. Notes in Math.,
Springer-Verlag,\textbf {1016}, (1983).

\bibitem{K6}
Kashiwara M., {\em D-modules and Microlocal Calculus,}
Translations of Mathematical Monographs, AMS, Vol \textbf{217},
(2003).

\bibitem{KO}
Kashiwara M. and Oshima T., {\em Systems of differential equations
with regular singularities and their boundary value problems,}
Ann. of Math.,\textbf{106},145-200,(1977).

\bibitem{KMS1}
Kashiwara M., Monteiro Fernandes T. and Schapira P., {\em
Truncated microsupport and holomorphic solutions of
$\mathcal{D}$-modules}, Ann. Scient. $\acute{E}c$. Norm. Sup.,
$4^e$ s$\acute{e}$rie, \textbf{36}, 583-599, (2003).

\bibitem{KMS2}
Kashiwara M., Monteiro Fernandes T. and Schapira P., {\em
Involutivity of truncated microsupports,} Bull. Soc. Math. France,
\textbf{131}(2), p.259-266, (2003)


\bibitem{K-S1}
Kashiwara M. and Schapira P., {\em Sheaves on manifolds,}
Grundlehren der Math. Wiss.,\textbf{292}, Springer Verlag (1990)

\bibitem{KS2}
Kashiwara M. and Schapira P., {\em Microlocal study of sheaves,}
Ast\' erisque \textbf{128} (1985).

\bibitem{LS}
Laurent Y.  and Schapira P., {\em Images inverses des modules
diff\' erentiels}, Compositio Math., \textbf{61},229-251,(1987).

\bibitem{S-K-K}
Sato M., Kawai T. and Kashiwara M.,  {\em Hyperfunctions and
pseudodifferential equations,} Lecture Notes in Math., Springer
\textbf{287}, 265-529 (1973).

\bibitem{S}
Schapira P., {\em Microdifferential systems in the complex domain}
Grundlehren der Math. Wiss., Vol. \textbf{269}, Springer-Verlag,
1985.

\bibitem{S2}
Schapira P., {\em Vanishing in Highest Degree for Solutions of
D-modules and Perverse Sheaves,} Publications of the Research
Institute for Mathematical Sciences, Kyoto univ.,
\textbf{26},3,535-538, (1990).

\bibitem{T}
Tonin F., {\em Holomorphic extension for solutions of the
characteristic Cauchy Problem} Preprint, (1998).

\end{thebibliography}
\end{document}